\documentclass[12pt,a4paper]{article}
\usepackage{hyperref}
\usepackage[german,english]{babel} 
\usepackage{amsthm}
\usepackage{amsmath}
\usepackage{amsfonts}
\usepackage{amssymb}
\usepackage{enumerate}
\usepackage{amscd}
\usepackage{graphicx}
 \usepackage{bbm}
\usepackage{makeidx}
\usepackage[all]{xy}
 
\newcounter{Def}[section]

\newtheorem{Examples}[Def]{Examples}
\newtheorem{Example}[Def]{Example}
\newtheorem{Remarks}[Def]{Remarks}

\newtheorem{Theorem}[Def] {Theorem}
\newtheorem{Kor}[Def] {Corollary}

\newenvironment{Proof}[1][\quad]
{\mbox{}\\[-17pt]\noindent \textbf{Proof:{\hspace{0.7cm}
#1}\\}}{\hfill$\Box$\\}

\newtheorem{Definition}[Def] {Definition}

\def\R          {\mathbb R}

\def\Z          {\mathbb Z}

\def\be                 {\begin{equation}}
\def\ee                 {\end{equation}}

\def\J          {\mathrm{J}^{\infty}\pi}

\def\id               {{\rm id}}

\def\id                 {\mathrm{id}}

\def\teins              {{\mathbf{1}}}
\def\c                  {\mathcal{C}}   
\def\d                  {\mathcal{D}} 
\def\h                  {\mathcal{H}}     
\def\hom                {\mathrm{Hom}}

\def\apo           {{\rm S}}
\def\Bimod         {\mbox{-bimod}}
\def\bP            {\begin{picture}}
\def\C             {{\ensuremath{\mathcal C}}}
\def\cir           {\,{\circ}\,}
\def\complex       {{\ensuremath{\mathbbm C}}}
\def\Enumerate     {\def\leftmargini{1.35em}~\\[-1.66em]\begin{enumerate}\addtolength\itemsep{-6pt}}
\def\eP            {\end{picture}}
\def\eps           {\varepsilon}
\def\eq            {\,{=}\,}
\def\findim        {fini\-te-di\-men\-si\-o\-nal}
\def\H             {\ensuremath{\mathcal H}}
\def\iN            {\,{\in}\,}
\newcommand\Includepichtft[1] {{\begin{picture}(0,0)(0,0)
                   \scalebox{.38}{\includegraphics{pic_htft_#1.eps}}\end{picture}}}
\def\ko            {\ensuremath{\Bbbk}}
\def\Mod           {\mbox{-mod}}
\def\one           {{\bf1}}
\def\oti           {\,{\otimes}\,}
\def\sse           {\scriptstyle}
\def\ssse          {\scriptscriptstyle}
\def\slz           {\ensuremath{\mathrm{SL}(2,\mathbb Z)}}
\def\Times         {\,{\times}\,}
\def\To            {\,{\to}\,}
\def\Vectcf        {\ensuremath{{\mathcal V}\!\mbox{\sl ect}_{\rm fin}(\complex)}}
\def\Vectk         {\ensuremath{{\mathcal V}\!\mbox{\sl ect}(\ko)}}
\def\Vectkf        {\ensuremath{{\mathcal V}\!\mbox{\sl ect}_{\rm fin}(\ko)}}

\begin{document}


\thispagestyle{empty}
\begin{flushright}
   {\sf ZMP-HH/10-8}\\
   {\sf Hamburger$\;$Beitr\"age$\;$zur$\;$Mathematik$\;$Nr.$\;$366}\\[2mm]
   March 2010
\end{flushright}
\vskip 2.0em
\begin{center}\Large 
HOPF ALGEBRAS AND FROBENIUS ALGEBRAS \\[3pt] IN FINITE TENSOR CATEGORIES
\end{center}\vskip 1.4em
\begin{center}
  Christoph Schweigert $^{a}$
  ~~and~~ J\"urgen Fuchs $^b$
\end{center}

\vskip 6mm

\begin{center}\it
  $^a$ Organisationseinheit Mathematik, \ Universit\"at Hamburg\\
  Bereich Algebra und Zahlentheorie\\
  Bundesstra\ss e 55, \ D\,--\,20\,146\, Hamburg \\[7pt] 
  $^b$ Teoretisk fysik, \ Karlstads Universitet
  \\Universitetsgatan 21, \ S\,--\,\,651\,88 Karlstad
\end{center}

\vskip 3em

\begin{abstract} \noindent
We discuss algebraic and representation theoretic structures in braided tensor 
categories \C\ which obey certain finiteness conditions. Much interesting structure
of such a category is encoded in a Hopf algebra \H\ in \C. In particular, the
Hopf algebra \H\ gives rise to representations of the modular group \slz\
on various morphism spaces. We also explain how every symmetric special Frobenius 
algebra in a semisimple modular category provides additional structure related to 
these representations.
\end{abstract}

\setcounter{footnote}{0} \def\thefootnote{\arabic{footnote}} 


\section{Braided finite tensor categories}

Algebra and representation theory in semisimple ribbon categories has been 
an active field over the last decade, having applications to quantum groups, 
low-dimensional topology and quantum field theory. More recently, partly in
connection with progress in the understanding of logarithmic conformal field
theories, there has been increased interest in tensor categories that 
are not semisimple any longer, but still obey certain finiteness conditions
\cite{EO}.

Owing to the work of various groups (for some recent results see e.g.\
\cite{GT,NT}), examples of such categories are by now rather explicitly
understood, at least as abelian categories. In this section we  describe a 
class of categories that has received particular attention.
This will allow us to define the structure of a semisimple modular tensor
category. To extend the notion of modular tensor category to the non-semisimple
case requires further categorical constructions involving Hopf algebras
and coends; these will be introduced in section 2. These constructions also
provide representations of the modular group \slz\ on certain morphism spaces.
In section 3 we show that symmetric special Frobenius algebras in semisimple
modular tensor categories give rise to structures related to such
\slz-representations.

\medskip

Let \ko\ be an algebraically closed field of characteristic zero and
\Vectkf\ the category of \findim\ \ko-vector spaces.

\begin{Definition}\mbox{} 
\\
A \emph{finite category} \C\ is an abelian category enriched over \Vectkf\
with the following additional properties:
\Enumerate
\item Every object has finite length.
\item Every object $X\iN\c$ has a projective cover $P(X)\iN\c$.
\item The set $I$ of isomorphism classes of simple objects is finite. 
\end{enumerate}
\end{Definition}

It can be shown that an abelian category is a finite category if and
only if it is equivalent to the category of (left, say) modules over
a finite-dimensional \ko-algebra.

We will be concerned with finite categories that have additional structure.
First, they are tensor categories, i.e., for our purposes, sovereign
monoidal categories:

\begin{Definition}\mbox{} \\
A \emph{tensor category} over a field \ko\ is a \ko-linear abelian monoidal 
category $\c$ with simple tensor unit $\teins$ and with a left and a right 
duality in the sense of {\rm \cite[Def.\,XIV.2.1]{Ka}}, such that the category
is sovereign, i.e.\ the two functors
  $$
  ?^\vee,\, {{}^\vee}?:~~ \c\to\c^{\mathrm{opp}}
  $$
that are induced by the left and right dualities coincide. 
\\
Thus for any object $V\iN\c$ there exists an object $V^\vee\eq {{}^\vee}V \iN\c$ 
together with morphisms
  $$
  b_V:~~ \teins\to V\oti V^\vee \quad\text{ and }\quad~
  d_V:~~ V^\vee\oti V \to \teins 
  $$
(right duality) and 
  $$
  \tilde b_V:~~ \teins\to V^\vee\oti V \quad\text{ and }\quad~
  \tilde d_V:~~ V\oti V^\vee \to \teins 
  $$
(left duality), obeying the relations
  $$ 
  (\id_V \oti d_V)\circ(b_V\oti \id_V)= \id_V \quad\text{ and }
  \quad
  (d_V\oti \id_{V^\vee})\circ (\id_{V^\vee}\oti b_V) = \id_{V^\vee_{}}
  $$
and analogous relations for the right duality, and the duality functors not only
coincide on objects, but also on morphisms, i.e.\
  $$
  (d_V\oti\id_{U^\vee_{}}) \cir (\id_{V^\vee_{}}\oti f\oti \id_{U^\vee_{}})
  \cir (\id_{V^\vee_{}} \oti b_U)
  = (\id_{U^\vee_{}} \oti \tilde d_V) \cir (\id_{U^\vee_{}}\oti f\oti \id_{V^\vee_{}})
  \cir (\tilde  b_U\oti\id_{V^\vee_{}}) 
  $$
for all morphisms $f\colon\,U\To V$. 
\end{Definition}

To give an example, the category of finite-dimensional left modules over any 
\findim\ complex Hopf algebra $H$ is a finite tensor category.
As a direct consequence of the definition, the tensor product functor 
$\otimes$ is exact in both arguments. We will impose on the dualities the
additional requirement that left and right duality lead to the same 
cyclic trace $\,\mathrm{tr}\colon\, \mathrm{End}(U) \To \mathrm{End}(\teins)$, 
and thus to a dimension $\,\mathrm{dim}(U) \eq \mathrm{tr}(\id_U)$.

The categories of our interest have in addition a braiding:

\begin{Definition} \mbox{} 
\\
A \emph{braiding} on a tensor category \C\ is a natural isomorphism 
  $$
  c:\quad \otimes \to \otimes^{\mathrm{opp}}
  $$
that is compatible with the tensor product, i.e.\ satisfies
  $$
  c_{U\otimes V,W} = (c_{U,W}\oti\id_V) \circ(\id_U\oti c_{V,W} )
  \quad \text{ and } \quad
  c_{U,V\otimes W} = (id_V \oti c_{U,W}) \circ (c_{U,V}\oti\id_W ) \, .
  $$
\end{Definition} 

We choose a set $\{U_i\}_{i\in I}$ of representatives for the
isomorphism classes of simple objects and take the tensor unit
to be the representative of its isomorphism class, writing $\one\eq U_0$.

We are now ready to formulate the notion of a modular tensor 
category. Our definition will, however, still be preliminary, as it has 
the disadvantage of being sensible only for semisimple categories.

\begin{Definition}\label{def.S}\mbox{} 
\\
A semisimple \emph{modular tensor category} is a semisimple finite braided 
tensor category such that the matrix $(S_{ij})_{i,j\in I}$ with entries
  $$
  S_{ij}:= \mathrm{tr} (c_{U_j,U_i} \cir c_{U_i,U_j})
  $$
is non-degenerate.
\end{Definition}

Two remarks are in order:

\begin{Remarks}\mbox{} \\[-1.6em]
\begin{enumerate}

\item 
The representation categories of several algebraic
structures give examples of semisimple modular tensor categories:
\\[-1.66em]\begin{enumerate}\addtolength\itemsep{-2pt}
\item 
Left modules over connected factorizable ribbon weak Hopf algebras 
with Haar integral over an algebraically closed field {\rm \cite{nitv}}.
\item 
Local sectors of a finite $\mu$-index net of von Neumann algebras
on $\R$, if the net is strongly additive and split {\rm \cite{KLM}}.
\item 
Representations of selfdual $C_2$-cofinite vertex algebras 
with an additional finiteness condition on the homogeneous
components and which have semisimple representation categories {\rm \cite{Hu}}.
\end{enumerate}

\item
By the results of Reshetikhin and Turaev {\rm \cite{RT,Tu}}, every 
\,\complex-linear semisimple modular tensor category $\c$ provides a 
three-dimensional topological field theory, i.e.\ a tensor functor
  $$
  \mathrm{tft}_\c:\quad \mathrm{cobord}_{3,2}^\c \to \Vectcf \,.
  $$
Here $\mathrm{cobord}_{3,2}^\c$ is a category of three-dimensional cobordisms 
with embedded ribbon graphs that are decorated by objects and morphisms of $\c$.
\\
There are also various results for the case of non-semisimple modular categories. 
We refer to {\rm \cite{henni,lyub6,V}} for the construction of three-manifold 
invariants, to {\rm \cite{lyub6}} for the construction of representations of 
mapping class groups, and to {\rm \cite{KElu}} for an attempt to unify these 
constructions in terms of a topological quantum field theory defined on a double 
category of manifolds with corners.
\end{enumerate}
\end{Remarks}


\section{Hopf algebras, coends and modular tensor categories}

Our goal is to study some algebraic and representation theoretic structures in 
tensor categories of the type introduced above. To simplify the exposition,
we suppose that we have replaced the tensor category $\c$ by an equivalent 
strict tensor category.

\begin{Definition}\label{Algebra} \mbox{}
\\
A (unital, associative) \emph{algebra} in a (strict) tensor category $\c$ is
a triple consisting of an object $A\iN\c$, a multiplication morphism 
$m \iN \hom(A \oti A,A)$ and a unit morphism $\eta \iN \hom(\teins, A)$, subject
to the relations
  $$
  m \circ (m \oti \id_A) = m \circ (\id_A \oti m) \qquad \text{ and }\qquad
  m \circ (\eta \oti \id_A) = \id_A = m \circ (\id_A \oti \eta ) \,.
  $$
which express associativity and unitality.
\\
Analogously, a \emph{coalgebra} in $\c$ is a triple consisting of an object 
$C$, a comultiplication $\Delta\colon\, C\To C\oti C$ and a counit 
$\eps\colon\, C\To \teins$ obeying coassociativity and counit conditions.
\end{Definition} 

Similarly one generalizes other basic notions of algebra to the categorical 
setting and introduces modules, bimodules, comodules, etc. (For a more 
complete exposition we refer to \cite{FRS}.)

To proceed we observe that the multiplication of an algebra $A$ endows both 
$A$ itself and $A\oti A$ with the structure of an $A$-bimodule.
Further, if the category $\c$ is braided, then the object $A\oti A$ can
be endowed with the structure of a unital associative algebra by taking
the morphisms $(m\oti m)\cir(\id_A\oti c_{A,A}\oti\id_A)$ as the product
and $\eta\oti\eta$ as the counit.

\begin{Definition} 
~\\
Let $\c$ be a tensor category and $A\iN\c$ an object which is endowed 
with both the structure $(A,m,\eta)$ of a unital associative algebra and 
the structure $(A,\Delta,\eps)$ of a counital coassociative coalgebra.
\Enumerate
\item 
$(A,m,\eta,\Delta,\eps)$ is called a \emph{Frobenius algebra} iff 
$\,\Delta\colon A\To A\oti A$ is a morphism of bimodules.
\item 
$(A,m,\eta,\Delta,\eps)$ is called a \emph{bialgebra} iff 
$\,\Delta\colon A\To A\oti A$ is a morphism of unital algebras. 
\item
A bialgebra with an antipode $\,\apo\colon A\To A$ (with properties 
analogous to the classical case) is called a \emph{Hopf algebra}.
\end{enumerate}
\end{Definition}

To construct concrete examples of such structures, we recall
a few notions from category theory.

\begin{Definition}\mbox{} 
\\
Let $\c$ and $\d$ be categories and $F\colon\,\c^{\mathrm{opp}}\Times
\c\To \d$ be a functor.
\Enumerate
\item
For $B$ an object of $\d$, a \emph{dinatural transformation} 
$\varphi\colon F\,{\Rightarrow}\, B$ is a family of morphisms
$\varphi_X\colon\, F(X,X)\To B\,$ for every object $X\iN\c$ such that the diagram
  $$\xymatrix @R+8pt{
  F(Y,X)\ \ar^{F(\id_Y,f)}[rr]\ar_{F(f,\id_X)}[d]&&\ F(Y,Y)\ar^{\varphi_Y^{}}[d] \\
  F(X,X)\ \ar^{\varphi_X^{}}[rr] && \, B
  } $$
commutes for all morphisms $X\,{\stackrel f\to}\, Y$ in $\c$.
\item 
A \emph{coend} for the functor $F$ is a dinatural transformation 
$\,\iota\colon F\,{\Rightarrow}\, A\,$ with the universal property that any 
dinatural transformation $\varphi\colon F\,{\Rightarrow}\, B$ uniquely factorizes:
  $$\xymatrix @R+4pt{
  F(Y,X)\ \ar^{F(\id_Y,f)}[rr]\ar_{F(f,\id_X)}[d]
  && F(Y,Y)\ar^{\iota_Y}[d]  \ar @/^2pc/^{\varphi_Y^{}}[ddr]& \\
  F(X,X)\ar @/_2pc/_{\varphi_X^{}}[drrr] \ar^{\iota_X}[rr] && \, A \ar@{-->}[dr] & \\
  &&& \, B
  } $$
\end{enumerate}
\end{Definition}

If the coend exists, it is unique up to unique isomorphism.
It is denoted by $\int^X\! F(X,X)$.
The universal property implies that a morphism with domain $\int^X\! F(X,X)$ 
can be specified by a dinatural family of morphisms $X^\vee\oti X\to B$ for 
each object $X\iN\c$. 

We are now ready to formulate the following result.

\begin{Theorem} {\rm \cite{L2}}~
\\
In a finite braided tensor category $\c$, the coend
  $$
  \mathcal{H} := \int^X\! X^\vee\oti X
  $$
of the functor
  $$
  \begin{array}{rll}
  F:\quad \c^{\mathrm{opp}}\Times \c &\!\!\to\!\!& \c \\[3pt]
  (U,V)&\!\!\mapsto \!\!& U^\vee\oti V 
  \end{array} $$
exists, and it has a natural structure of a Hopf algebra in $\c$.
\end{Theorem}

\begin{Proof}
For a proof we refer e.g.\ to \cite{V}. Here we only indicate how the 
structural morphisms of the Hopf algebra are constructed. Owing to the universal
property, the counit $\eps_\h\colon\, \h\To\teins$ can be
specified by the dinatural family 
  $$ 
  \eps_\h \circ \iota_X = d_X :~~ X^\vee\oti X\to\teins
  $$
of morphisms. Similarly, the coproduct is given by the dinatural family
  $$
  \Delta_\h\circ\iota_\h = (\iota_X\oti\iota_X)\circ
  (\id_{X^\vee}\oti b_X\oti \id_X ) :~~ X^\vee\oti X\to\H\oti\H \,.
  $$
It should be appreciated that the braiding does not enter in
the coalgebra structure of \H. 
\\
It does enter in the product, though. We refrain from writing out the product 
as a formula. Instead, we use the graphical formalism \cite{joSt5,FRS} to 
display all structural morphisms $(m_\H,\Delta_\H,\eta_\H,\epsilon_\H,\apo_\H)$ 
of the Hopf algebra \H. More precisely, we display dinatural
families of morphisms so that the identities apply to all $X,Y\iN\C$:
  $$
  \bP(390,220) { \put(13,103) {
    \put(0,0) {
  \put(0,0)   {\Includepichtft{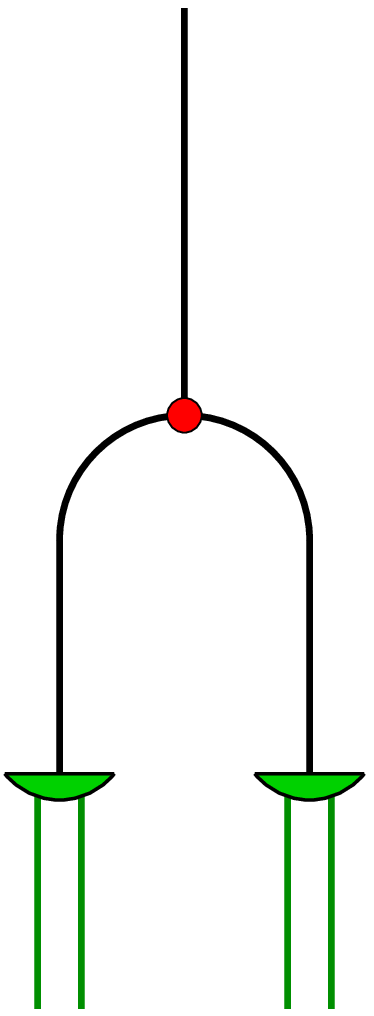}
  \put(-6,-9)   {$\sse X^{\!\vee} $}
  \put(7,-9)    {$\sse X $}
  \put(16.3,115){$\sse \H $}
  \put(22,69)   {$\sse m_\H $}
  \put(23,-9)   {$\sse Y^{\!\vee} $}
  \put(36,-9)   {$\sse Y $}
  \put(-9.6,23) {$\sse \iota_{\!X}^{} $}
  \put(41,23)   {$\sse \iota_{\!Y}^{} $}
  \put(60,50)   {$ = $}
   }
  \put(82,0) { {\Includepichtft{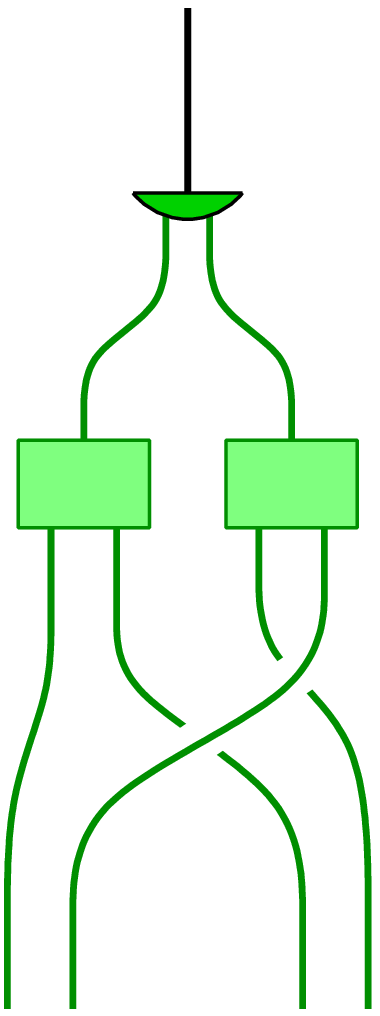}}
  \put(-10.8,66)  {$\sse \gamma_{X,Y} $}
  \put(34,65.2)   {$\sse \id_{Y|X} $}
  \put(-7,-8)   {$\sse X^{\!\vee} $}
  \put(6,-8)    {$\sse X $}
  \put(17.6,115){$\sse \H $}
  \put(25,-9)   {$\sse Y^{\!\vee} $}
  \put(38,-8)   {$\sse Y $}
  \put(-11,79)  {$\ssse (Y{\otimes}X)^{\!\vee}_{} $}
  \put(25.5,79)   {$\ssse Y{\otimes}X $}
   } }
    \put(210,8) {
  \put(0,0) { {\Includepichtft{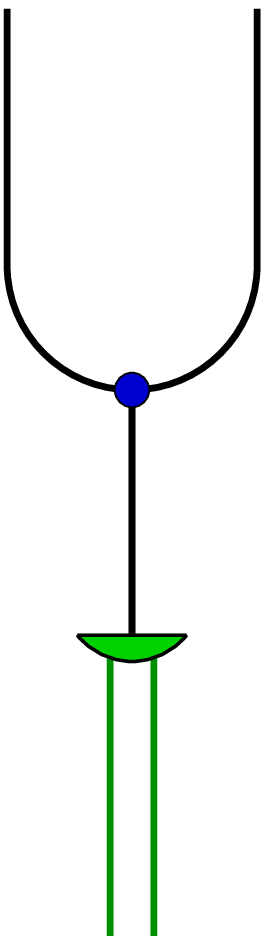}}
  \put(-3.4,107){$\sse \H $}
  \put(24.4,107){$\sse \H $}
  \put(2,-9)    {$\sse X^{\!\vee} $}
  \put(15,-9)   {$\sse X $}
  \put(15.9,53.5) {$\sse \Delta_\H $}
  \put(55,50)   {$ = $}
  \put(85,0) { {\Includepichtft{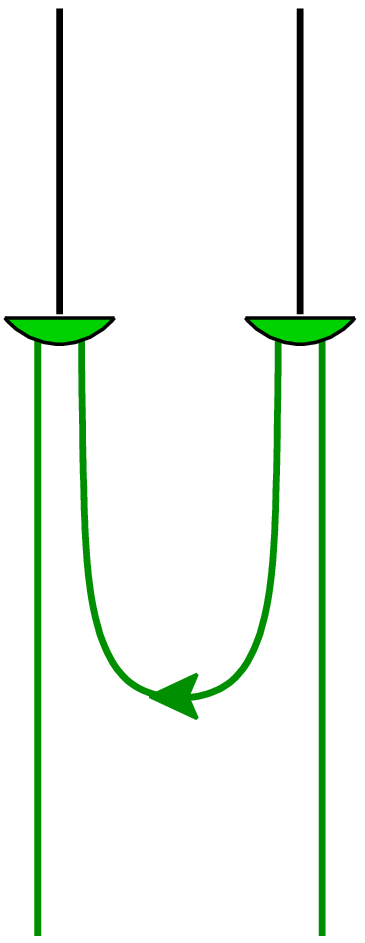}}
  \put(1.9,107) {$\sse \H $}
  \put(28.8,107){$\sse \H $}
  \put(-2,-9)   {$\sse X^{\!\vee} $}
  \put(30.5,-9) {$\sse X $}
   } } }
   }
    \put(13,0) {
  \put(0,17) { {\Includepichtft{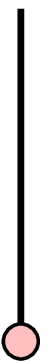}}
  \put(-1.8,42) {$\sse \H $}
  \put(5.7,1)   {$\sse \eta_\H^{} $}
  \put(25,23)   {$ = $}
  \put(52,-18) {\Includepichtft{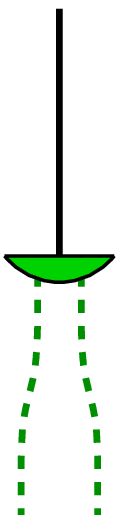}}
  \put(54,42)   {$\sse \H $}
   }
  \put(114,4) { {\Includepichtft{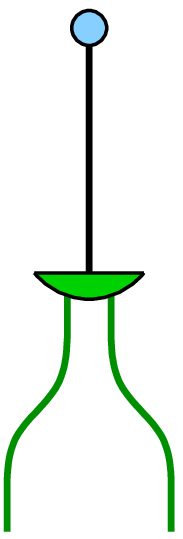}}
  \put(13.5,54.5) {$\sse \eps_\H^{} $}
  \put(-5.2,-9) {$\sse X^{\!\vee} $}
  \put(14,-9)   {$\sse X $}
  \put(38,30)   {$ = $}
  \put(67,10)   {\Includepichtft{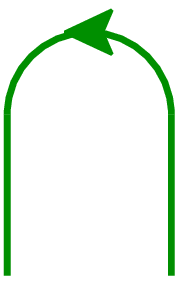}}
  \put(62,1)    {$\sse X^{\!\vee} $}
  \put(81,1)    {$\sse X $}
   }
  \put(260,-11) { {\Includepichtft{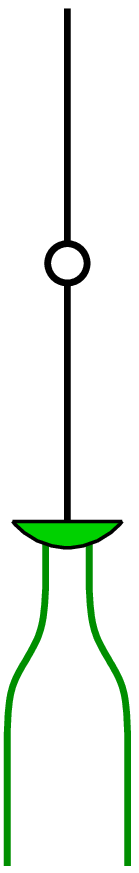}}
  \put(11.2,64.5) {$\sse \apo_\H $}
  \put(3.8,98)  {$\sse \H $}
  \put(34,44)   {$ = $}
  \put(63,0)   {\Includepichtft{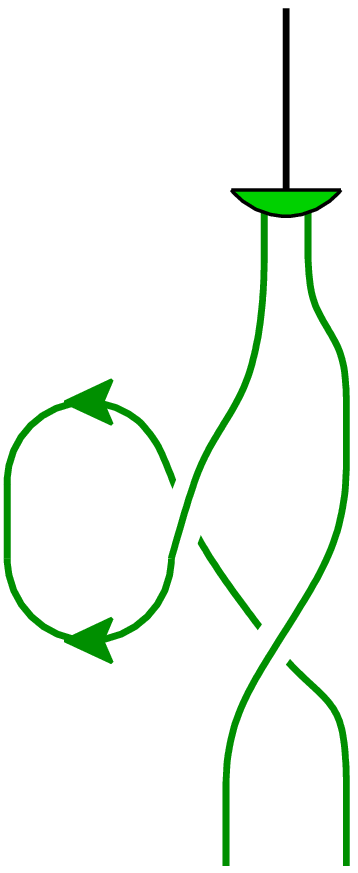}}
  \put(90.2,98) {$\sse \H $}
  \put(74.7,62) {$\ssse X^{\!\vee\!\vee} $}
  \put(98.9,62) {$\ssse X^{\!\vee} $}
   } } }
\eP
  $$
(Here $\gamma_{X,Y}$ is the canonical identification of $X^\vee{\otimes}\,Y^\vee$ 
with $(Y{\otimes}X)^\vee$, and $\id_{X|Y}$ is the one of $\id_X\oti \id_Y$ with 
$\id_{X\otimes Y}$.)
\end{Proof}

An explicit description of the Hopf algebra
$\h\iN\c$\ is available in the following specific situations:

\begin{Examples}
\Enumerate
\item 
For $\c\eq H\Mod$ the category of left modules over a finite-dimensional 
ribbon Hopf algebra $H$, the coend $\h \,\eq \int^X\! X^\vee\oti X$ is the 
dual space $H^*\eq\hom_\ko(H,\ko)$ endowed with the coadjoint representation.
The structure morphism for the coend for a module $M\iN H\Mod$ is
  $$
  \begin{array}{rll}
  \iota_M:\quad M^\vee\oti M &\!\! \to \!\!&H^* \\[3pt]
  \tilde m\oti m&\!\!\mapsto\!\!& \left(
  h\,{\mapsto}\, \langle\tilde m,h.m\rangle\right) \,.
  \end{array}
  $$
For more details see {\rm \cite[Sect.\,4.5]{V}}.
\item
If the finite tensor category \,\C\ is semisimple, then the Hopf algebra 
decomposes as an object as
$\,\h \eq \bigoplus_{i\in I} U_i^\vee\oti U_i $, see {\rm \cite[Sect.\,3.2]{V}}.
\end{enumerate}
\end{Examples}

The Hopf algebra in question has additional structure: it comes with an
integral and with a Hopf pairing.

\begin{Definition}\mbox{} 
\\
A \emph{left integral} of a bialgebra $(H,m,\eta,\Delta,\eps)$ in \C\ is a 
non-zero morphism $\mu_l\iN\hom(\teins,H)$ satisfying
  $$
  m\circ (\id_H\oti\mu_l) = \mu_l\circ \eps\,. 
  $$
A \emph{right cointegral} of $H$ is a non-zero morphism $\lambda_r\iN\hom(H,\one)$
satisfying
  $$
  (\lambda\oti\id_H) \circ \Delta =  \eta \circ \lambda \,.
  $$
Right integrals $\mu_r$ and left cointegrals $\lambda_l$ are defined analogously.
\end{Definition}

The Hopf algebra $\h$ in any finite braided tensor category has left and
right integrals, as can be shown \cite{L2} by a generalization of the 
classical argument of Sweedler that an integral exists for any finite-dimensional 
Hopf algebra. If \,\C\ is semisimple, then the integral of \H\ can be given 
explicitly \cite[Sect.\,2.5]{kerl5}:
  $$
  \mu_l = \mu_r = \bigoplus_{i\in I} \dim(U_i) \, b_{U_i} \,. 
  $$

\begin{Remarks} \label{2.7}
\Enumerate
\item
If the left and right integrals of \H\ coincide, then the integral can be used 
as a Kirby element and provides invariants of three-manifolds {\rm \cite{V}}. If
the category $\c$ is the category of representations of a finite-dimensional 
Hopf algebra, this is the Hennings-Lyubashenko {\rm \cite{lyub6}} invariant.

\item
The category $\c$ is semisimple if and only if the morphism
$\eps\cir\mu\iN \hom(\teins,\teins)$ does not vanish, 
i.e.\ iff the constant $\mathcal D^2$ of proportionality in
  $$
  \eps\circ\mu = \mathcal D^2 \, \id_{\teins}
  $$ 
is non-zero. (This generalizes Maschke's theorem.) This constant, in turn,
which in the semisimple case (with $\mu_l\eq\mu_r$ normalized as above) 
has the value $\mathcal D^2 \eq \sum_{i\in I}(\dim U_i)^2$,
crucially enters the normalizations in the Reshetikhin-Turaev construction
of topological field theories (see e.g.\ chapter II of {\rm \cite{Tu}}).
\\
Invariants based on nonsemisimple categories, like the Hennings invariant, 
vanish on many three-manifolds. This can be traced back to the vanishing of
$\eps\cir\mu$ {\rm \cite{CKS}}.

\item
Any Hopf algebra $H$ in $\c$ with invertible antipode that has a left integral 
$\mu$ and a right cointegral $\lambda$ with $\lambda\cir\mu \,{\ne}\, 0$
is naturally also a Frobenius algebra, with the same algebra structure. 

\end{enumerate}
\end{Remarks}

\begin{Definition}\mbox{} 
\\
A \emph{Hopf pairing} of a Hopf algebra $H$ in \C\ is a morphism
  $$
  \omega_H:\quad H\oti H \to \teins
  $$
such that
  $$ 
  \begin{array}l
  \omega_H\circ (m\oti \id_H)
  = (\omega_H\oti\omega_H)\circ (\id_H\oti c_{H,H}^{}\oti\id_H) \circ
  (\id_H\oti\id_H\oti\Delta) \,,
  \\[6pt]
  \omega_H\circ (\id_H\oti m)
  = (\omega_H\oti\omega_H)\circ (\id_H\oti c_{H,H}^{-1}\oti\id_H) \circ
  (\Delta\oti\id_H\oti\id_H)
  \\[6pt]
  \mbox{and}\qquad
  \omega_H\circ (\eta\oti \id_H) = \eps = \omega_H\circ (\id_H\oti\eta) \,.
  \end{array}
  $$
\end{Definition}

As one easily checks, a non-degenerate Hopf pairing gives an isomorphism
$H\To H^\vee$ of Hopf algebras. 

The dinatural family of morphisms
  $$
  (d_X \oti d_Y) \circ
  [\id_{X^\vee} \oti (c_{Y^\vee,X}^{}\cir c_{X,Y^\vee}^{}) \oti\id_Y]
  $$
induces a bilinear pairing $\omega_\H^{}\colon \H\oti\H\To\one$ on the
coend $\H \eq \int^X\! X^\vee\oti X$ of a finite braided tensor category.
It endows \cite{lyub6} the Hopf algebra \H\  with a symmetric Hopf pairing.

\medskip

We are now finally in a position to give a conceptual definition
of a modular finite tensor category without requiring it
to be semisimple:

\begin{Definition}{\rm \cite[Def.\ 5.2.7]{KElu}} \label{nssi-mtc}\mbox{} \\
A \emph{modular finite tensor category} is a braided finite tensor category 
for which the Hopf pairing $\omega_\h$ is non-degenerate.
\end{Definition}

\begin{Example}\mbox{} \\
The category $H\Mod$ of left modules over a finite-dimensional factorizable
ribbon Hopf algebra $H$ is a modular finite tensor category {\rm \cite{lyma,lyub6}}.
\end{Example}

One can show \cite[Thm.\,6.11]{L2} that if $\c$ is modular in the sense of definition
\ref{nssi-mtc}, then the left integral and the right integral of $\h$ coincide.

\medskip

As the terminology suggests, there is a relation with the modular group
\slz. To see this, we will now obtain elements $S_\h,T_\h\iN \mathrm{End}(\H)$
that satisfy the relations for generators of \slz.

Recall the notion of the center $Z(\C)$ of a category as the algebra of natural 
endotransformations of the identity endofunctor of \C\ \cite{maji16}. Given such 
a natural transformation $(\phi_X)_{X\in\C}$ with $\phi_X\iN\mathrm{End}(X)$, one 
checks that $(\iota_X\cir(\id_{X^\vee_{}}\oti\phi_X))_{X\in\C}$ is a dinatural 
family, so that the universal property of the coend gives us a unique 
endomorphism $\overline\phi_\H$ of \H\ such that the diagram
  $$\xymatrix @R+4pt{
  X^\vee\oti X\ \ar^{\id\otimes\phi_X^{}}[rr]\ar_{\iota_X}[d]
  && \, X^\vee\oti X \ar^{\iota_X}[d]\\
  \h\, \ar@{-->}_{\overline\phi_\H^{}}[rr]&&\, \h
  }$$
commutes, leading to an injective linear map $Z(\C) \To \mathrm{End}(\H)$. 

Since \H\ has in particular the structure of a coalgebra and $\one$ the
structure of an algebra, the vector space $\hom(\H,\one)$ has a natural 
structure of a \ko-algebra. Concatenating with the counit $\eps_\H$ gives a map
  $$
  Z(\C) \longrightarrow \mathrm{End}(\H) \stackrel{{(\eps_\H^{})}_*}
  {-\!\!\!-\!\!\!\longrightarrow} \hom(\H,\one) \,, 
  $$
which can be shown \cite[Lemma\,4]{kerl5} to be an isomorphism
of \ko-algebras. The vector space on the right hand side is dual to the 
vector space $\hom(\one,\H)$, of which one can think as the appropriate 
substitute for the space of class functions. Hence $\hom(\one,\H)$ would be a natural 
starting point for constructing a vector space assigned
to the torus $\mathrm T^2$ by a topological field theory based on \C.

\medskip

If the category $\c$ is a ribbon category, we have the ribbon element 
$\nu\iN Z(\C)$. We set
  $$
  T_\H := \overline\nu_\H^{} \,\,\in\,\mathrm{End}(\h) \,.
  $$
Pictorially,
  $$
  \bP(96,75)(0,-5) {
  \put(0,0)   {\Includepichtft{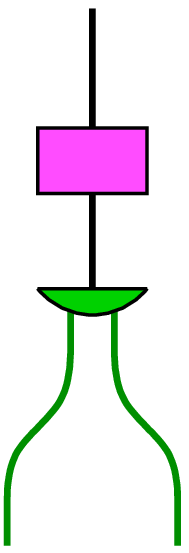}
  \put(-11.1,40.7) {$\overline\nu_\H^{}$}
  \put(-6,-9)   {$\sse X^{\!\vee} $}
  \put(14.9,-9) {$\sse X $}
  \put(6.8,62.7){$\sse \H $}
  }
  \put(46,28) {$=$}
  \put(77,0)  {\Includepichtft{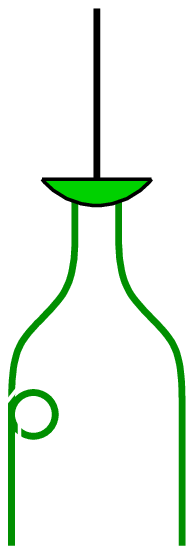}
  \put(-3,-9)   {$\sse X^{\!\vee} $}
  \put(17.9,-9) {$\sse X $}
  \put(8.4,62.7){$\sse \H $}
  } }
  \eP
  $$

Another morphism $\varSigma\colon \H\oti\H\To\H$ is obtained from the following 
family of morphisms which is dinatural both in $X$ and in $Y$:
  $$
  \bP(154,86)(0,-2) {
  \put(0,0)   {\Includepichtft{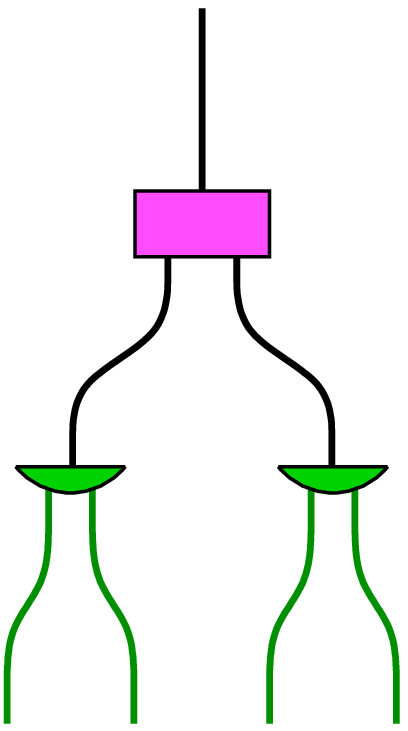}
  \put(3.5,51.6) {$\varSigma$}
  \put(-6,-9)   {$\sse X^{\!\vee} $}
  \put(9.9,-9)  {$\sse X $}
  \put(25,-9)   {$\sse Y^{\!\vee} $}
  \put(40.2,-9) {$\sse Y $}
  \put(18.8,83) {$\sse \H $}
  }
  \put(69,35) {$:=$}
  \put(108,0) {\Includepichtft{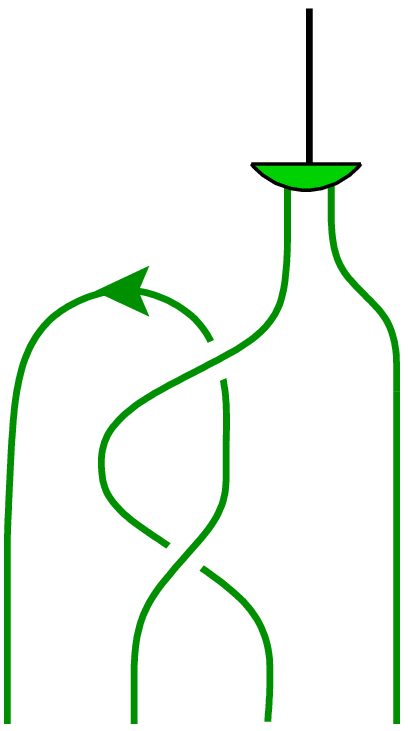}
  \put(-6,-9)   {$\sse X^{\!\vee} $}
  \put(9.9,-9)  {$\sse X $}
  \put(25,-9)   {$\sse Y^{\!\vee} $}
  \put(40.2,-9) {$\sse Y $}
  \put(29.8,83) {$\sse \H $}
  } }
  \eP
  $$
Composing this morphism to \H\ with a left or right integral 
$\mu\colon \one\To\H$ one arrives at an endomorphism 
  $$
  S_\H := \varSigma \circ (\id_\H \oti \mu) \,\,\iN \,\mathrm{End}(\H) \,.
  $$

For $\xi\iN \ko^\times$, denote by $\ko_\xi\slz$ the twisted group algebra 
of \slz\ with relations $S^4\eq 1$ and $(ST)^3\eq \xi\, S^2$.
The previous construction and the following result are due to Lyubashenko.

\begin{Theorem} {\rm \cite[Sect.\,6]{L2}}~
\\
Let \C\ be modular. Then the two-sided integral of \H\ can be normalized 
in such a way that the endomorphisms $S_\H$ and $T_\H$ of \H\ provide
a morphism of algebras
  $$
  \ko_\xi\slz \longrightarrow \mathrm{End}(\H)
  $$
for some $\xi\iN \ko^\times$.
\end{Theorem}

Since for every $U\iN\C$ the morphism space $\hom(U,\H)$ is, by push-forward, 
a left module over the algebra $\mathrm{End}(\H)$, we obtain this way
projective representations of \slz\ on all vector spaces $\hom(U,\H)$.

\medskip

To set the stage for the results in the next section, we consider the map
  $$
  \begin{array}{rll}
  \mathrm{Obj}(\c) &\!\!\to \!\!& \hom(\teins,\h) \\[2pt]
  U&\!\!\mapsto \!\!& \chi_U^{} 
  \end{array}
  $$
with 
  $$
  \chi_U^{}:\quad \teins \stackrel{b_U}\longrightarrow U^\vee\oti U 
  \stackrel{\iota_U^{}}\longrightarrow \h \,.
  $$
It factorizes to a morphism of rings
  $$
  K_0(\c)\to \hom(\teins,\h) = \mathrm{tft}_{\c}(T^2)\,.
  $$
If the category $\c$ is semisimple, then
$\hom(\teins,\h) \,{\cong}\, \bigoplus_{i\in I} \hom(\teins,U^\vee_i\oti U_i)$,
so that $\{\chi_{U_i}^{}\}_{i\in I}^{}$ constitutes a basis of the vector 
space $\hom(\teins,\h)$. If $\c$ is not semisimple, these elements are still 
linearly independent, but they do not form a basis any more. Pseudo-characters
\cite{Mi,GT} have been proposed as a (non-canonical)
complement of this linearly independent set.


\section{Frobenius algebras and braided induction}

In this section we show that symmetric special Frobenius algebras (i.e.\
Frobenius algebras with two further properties, to be defined below) in a 
modular tensor category allow one to specify interesting
structure related to the \slz-representation that we have just explained.

Given an algebra $A$ in a braided (strict) tensor category, we consider 
the two tensor functors
  $$
  \begin{array}{rll}
  \alpha_{\!A}^\pm:\quad \c&\!\! \to\!\!& A\Bimod \\[2pt]
  U&\!\!\mapsto\!\!& \alpha^\pm_{\!A}(U)\ 
  \end{array}
  $$
which assign to an object $U\in\c$ the bimodule $(A\oti U,\rho_l,\rho_r)$ 
for which the left action is given by multiplication and the right action by
multiplication composed with a braiding,
  $$ 
  \begin{array}{ll}&
  \rho_l = m\oti \id_U\in \hom(A\oti A\oti U,A\oti U) \\[7pt]{\rm and}\quad&
  \rho^+_r = (m\oti \id_U)\circ (\id_A\oti c_{U,A}) \,\, \text{ and } \,\,
  \rho^-_r = (m\oti \id_U)\circ (\id_A\oti c_{A,U}^{-1}) 
\,.
  \end{array}
  $$
We call these functors \emph{braided induction} functors. They have been 
introduced, under the name $\alpha$-induction, in operator algebra theory 
\cite{LR,X,BE}. For more details in a category-theoretic framework we refer 
to \cite[Sect.\,5.1]{O}.

We pause to recall that \cite{VZ} an \emph{Azumaya algebra} $A$ is an 
algebra for which the two functors $\alpha_{\!A}^\pm$ are equivalences 
of tensor categories. This should be compared to the textbook definition of
an Azumaya algebra in the tensor category of modules over a commutative 
\ko-algebra $A$, requiring in particular the morphism
  $$
  \begin{array}{rll}
  \psi_A:\quad A\oti A^{\rm opp} &\!\!\to \!\!&\mathrm{End}(A) \\[2pt]
  a\otimes a'&\!\!\mapsto \!\!& \left(x\mapsto a\,{\cdot}\, x\,{\cdot}\, a'\right)
  \end{array}
  $$ 
to be an isomorphism of algebras. Indeed, in this situation for an Azumaya
algebra $A$ one has the following chain of equivalences:
  $$
  A\Bimod \stackrel\sim\longrightarrow A{\otimes}A^{\rm opp}\Mod
  \stackrel{\psi_A}{-\!\!\!\longrightarrow} \mathrm{End}(A)\Mod
  \stackrel{\mathrm{Morita}}{-\!\!\!-\!\!\!-\!\!\!-\!\!\!\longrightarrow} \Vectk \,.
  $$

We now introduce the properties of an algebra $A$ to be symmetric and special.

\begin{Definition}\label{Frobenius}\mbox{} 
\\
Let $\c$ be a tensor category. 
\Enumerate 

\item 
For \C\ enriched over the category of \ko-vector spaces, a 
\emph{special algebra} in $\c$ is an object $A$ of \C\ that is endowed
with an algebra structure $(A,m,\eta)$ and a coalgebra structure 
$(A,\Delta,\eps)$ such that 
  $$
  \eps \circ \eta  = \beta_1\, \id_{\teins}
  \quad\text{ and }\quad m \circ \Delta  = \beta_{A}\, \id_{A}
  $$
with invertible elements $\beta_1, \beta_A \iN \ko^{\times}$.
    
\item 
A \emph{symmetric algebra} in $\c$ is an algebra $(A,m,\eta)$ together with
a morphism $\eps \iN \hom(A,\teins)$ such that the two morphisms 
  \begin{align}
  \Phi_1 & := [(\eps \circ m) \oti \id_{A^\vee_{}}] \circ (\id_A \oti b_A)
  \,\in\hom(A,A^{\vee})
  \qquad{\rm and} \\[.2em]
  \Phi_2 & := [\id_{A^\vee_{}} \oti (\eps \circ m)] \circ (\tilde{b}_A \oti \id_A)
  \,\in\hom(A,A^{\vee})
  \end{align}
are identical.
\end{enumerate}
\end{Definition}

Special algebras are in particular separable, and as a consequence their 
categories of modules and bimodules are semisimple. A class of examples of 
special Frobenius algebras is supplied by the Frobenius algebra structure on
a Hopf algebra $H$ in $\c$, provided $H$ is semisimple.

\medskip

We now consider the case of a semisimple modular tensor category \C\ and introduce 
for any algebra $A$ in $\c$ the  square matrix $(Z_{ij})_{i,j\in I}$ with entries
  $$
  Z_{ij}(A):= \dim_\ko\hom_{A|A}^{}(\alpha_{\!A}^-(U_i^{}),\alpha_{\!A}^+(U_j^\vee))\,,
  $$
where $\hom_{A|A}$ stands for homomorphisms of bimodules. Identifying
$A\Bimod$ with the tensor category of module endofunctors of $A\Mod$, one sees
that the non-negative integers $Z_{ij}(A)$ only depend on the Morita class of $A$.

In this setting, and in case that the algebra $A$ is symmetric and special,
we can make the following statements.

\begin{Theorem} {\rm \cite[Thm.\,5.1(i)]{FRS} }~ \\
For \,\C\ a semisimple modular tensor category and $A$ a special symmetric
Frobenius algebra in \C, the morphism
  \begin{equation}
  \sum_{i,j\in I}Z_{ij}(A) \, \chi_i^{}\oti\chi_j^{} 
  \,\in \hom(\teins,\h) \,{\otimes_\ko}\, \hom(\teins,\h)
  \label{eqZ}\end{equation}
is invariant under the diagonal action of \,\slz.
\end{Theorem}

\begin{Remarks}\mbox{} \\[-1.6em]
\begin{enumerate}
\item
In conformal field theory, the expression (\ref{eqZ}) has the interpretation
of a partition function for bulk fields.
\\[-1.7em]~

\item 
For semisimple tensor categories based on the $\mathrm{sl}(2)$ affine
Lie algebra, an A-D-E pattern appears {\rm \cite{kios}}.
\end{enumerate}
\end{Remarks}

We finally summarize a few other results that hold under the assumption
that $\c$ is a semisimple modular tensor category and $A$ a symmetric 
special Frobenius algebra in $\c$.  To formulate them, we need the following 
ingredients: The \emph{fusion algebra} 
  $$
  R_\C := K_0(\c) \,{\otimes_\Z}\, \ko
  $$
is a separable commutative algebra with a natural basis $\{\,[U_i]\,\}_{i\in I}$ 
given by the isomorphism classes of simple objects. The matrix $S$ introduced 
in definition \ref{def.S} provides a natural bijection from the set
of isomorphism classes of irreducible representations of $R_\C$ to $I$.

\begin{Theorem} {\rm \cite[Thm.\,5.18]{FRS}}~ 
\\
For any special symmetric Frobenius algebra $A$ the vector space 
$K_0(A\Mod)\,{\otimes_\Z}\,\ko$ is an $R_\C$-module. The multiset 
$\,\mathrm{Exp}(A\Mod)$ that contains the irreducible $R_\C$-representations, 
with their multiplicities in this $R_\C$-module, can be expressed in terms of 
the matrix $Z(A)$:
  $$
 \mathrm{Exp}(A\Mod) \,=\, \mathrm{Exp}(Z(A))
 := \{i\iN I\ \text{with multiplicity } Z_{ii}(A)\} \,.
  $$ 
\end{Theorem}

The observation that the vector space $K_0(A\Mod)\,{\otimes_\Z}\,\ko$
has a natural basis provided by the classes of simple $A$-modules gives 

\begin{Kor}\mbox{} \\
The number of isomorphism classes of simple $A$-modules equals $\,\mathrm{tr}(Z(A))$.
\end{Kor}

The category $A\Bimod$ of $A$-bimodules has the structure of a tensor 
category. From the fact that $A$ is a symmetric special Frobenius algebra, 
it follows \cite{FS} that $A\Bimod$ inherits left and right dualities 
from $\c$. Hence the tensor product on $A\Bimod$ is exact and thus 
$K_0(A\Bimod)$\ is a ring. The corresponding \ko-algebra can again
be described in terms of the matrix $Z(A)$:

\begin{Theorem} {\rm \cite{O,fuRs12}}~
\\
There is an isomorphism 
  $$
  K_0(A\Bimod) \,{\otimes_\Z}\, \ko \,\cong\,
  \bigoplus_{i,j\in I} {\rm Mat}_{Z_{ij}(A)}(\ko)\,,
  $$
of \ko-algebras, with ${\rm Mat}_n(\ko)$ denoting the algebra of 
\,\ko-valued $n\Times n$-matrices.
\end{Theorem}

\begin{Kor}\mbox{} \\
The number of isomorphism classes of simple $A$-bimodules
equals $\,\mathrm{tr}( ZZ^{\rm t})$.
\end{Kor}

\begin{Theorem} {\rm \cite[Prop.\,4.7]{ffrs5}}~
\\
Any $A$-bimodule is a subquotient of a bimodule of the form
$\alpha_{\!A}^+(U) \,{\otimes_{\!A}}\, \alpha_{\!A}^-(V)$ for some pair
of objects $U,V\iN\c$.
\end{Theorem}

\section{Outlook}

We conclude this brief review with a few comments. First, all the results 
about algebra and representation theory in braided tensor categories 
that we have presented above are motivated by a construction of correlation
functions of a rational conformal field theory as elements of vector spaces 
which are assigned by a topological field theory to a two-manifold. For 
details of this construction we refer to \cite{SFR} and the literature
given there.

In the conformal field theory context the matrix $Z$ describes the partition
function of bulk fields. The three-dimensional topology involved in the RCFT 
construction provides in particular a motivation for using the different 
braidings which lead to the functors $\alpha_{\!A}^+$ and $\alpha_{\!A}^-$ as 
well as in the definition of $Z(A)$.

To extend the results obtained in connection with rational conformal field 
theory to non-semisimple finite braided tensor categories remains a major 
challenge. Intriguing first results include, at the level of chiral data, a 
generalization of the Verlinde formula (see \cite{GT} and references given
there), and at the level of partition functions, the bulk partition functions 
for logarithmic conformal field theories in the $(1,p)$-series found in \cite{GR}.

\vskip3em

\noindent{\sc Acknowledgments:}\\[2pt]
JF is partially supported by VR under project no.\ 621-2009-3343. 
\\
CS is partially supported by the DFG Priority Program 1388 ``Representation theory''.

\newpage

 \newcommand\wb{\,\linebreak[0]} \def\wB {$\,$\wb}
 \newcommand\Bi[2]    {\bibitem[#2]{#1}} 
 \newcommand\inBO[9]  {{\em #9}, in:\ {\em #1}, {#2}\ ({#3}, {#4} {#5}), p.\ {#6--#7} {{\tt [#8]}}}
 \newcommand\JO[6]    {{\em #6}, {#1} {#2} ({#3}) {#4--#5} }
 \renewcommand\J[7]   {{\em #7}, {#1} {#2} ({#3}) {#4--#5} {{\tt [#6]}}}
 \newcommand\K[6]     {{\em #6}, {#1} {#2} ({#3}) {#4} {{\tt [#5]}}}
 \newcommand\BOOK[4]  {{\em #1\/} ({#2}, {#3} {#4})}
 \newcommand\prep[2]  {{\em #2}, preprint {\tt #1}}
 \def\adma  {Adv.\wb Math.}
 \def\apcs  {Applied\wB Cate\-go\-rical\wB Struc\-tures}
 \def\comp  {Com\-mun.\wb Math.\wb Phys.}
 \def\duke  {Duke\wB Math.\wb J.}
 \def\fiic  {Fields\wB Institute\wB Commun.}
 \def\ijmp  {Int.\wb J.\wb Mod.\wb Phys.\ A}
 \def\inma  {Invent.\wb math.}
 \def\joal  {J.\wB Al\-ge\-bra}
 \def\jopa  {J.\wb Phys.\ A}
 \def\jlms  {J.\wB London\wB Math.\wb Soc.}
 \def\jpaa  {J.\wB Pure\wB Appl.\wb Alg.}
 \def\momj  {Mos\-cow\wB Math.\wb J.} 
 \def\mpcp  {Math.\wb Proc.\wB Cam\-bridge\wB Philos.\wb Soc.}
 \def\nupb  {Nucl.\wb Phys.\ B}
 \def\plms  {Proc.\wB Lon\-don\wB Math.\wb Soc.}
 \def\pnas  {Proc.\wb Natl.\wb Acad.\wb Sci.\wb USA} 
 \def\rvmp  {Rev.\wb Math.\wb Phys.}
 \def\taia  {Top\-o\-lo\-gy\wB and\wB its\wB Appl.}
 \def\trgr  {Trans\-form.\wB Groups}

\end{document}